\documentclass[envcountsame]{llncs}
\usepackage[utf8]{inputenc}
\usepackage{amssymb}
\usepackage{pstricks,pst-node,pstcol}

\newcommand{\conju}[2]{\mathcal{#1}^{\mathcal{F}}(#2)}
\newcommand{\bru}{B^{\mathcal{F},n}}
\newcommand{\gbru}[1]{G^{\mathcal{F},#1}}
\newcommand{\Zset}{\bbbz}
\newcommand{\mod}{\mathrm{\, mod\,}}
\spnewtheorem*{procedure}{Procedure}{\bfseries}{\itshape}

\begin{document}

\title{Minimum de Bruijn Sequence in a Language with Forbidden
  Substrings%
\thanks{Partially supported by ECOS C00E03 (French-Chilean Cooperation), 
Programa Iniciativa Cient\'{\i}fica Milenio P01-005,
and CONICYT Ph.D. Fellowship.}
}
\author{Eduardo Moreno\inst{1} \and Mart\'{\i}n
  Matamala\inst{2}\inst{3}}
\institute{%
Universidad Adolfo Ibañez\\ Avda.  Diagonal las Torres 2640, Peñalolen, Santiago, Chile.\\
\email{eduardo.moreno@uai.cl}
\and
Departamento de Ingenier\'{\i}a Ma\-te\-m\'a\-ti\-ca, Facultad de
Ciencias F\'{\i}sicas y Ma\-te\-m\'a\-ti\-cas, Universidad de Chile. Casilla
170-3, Correo 3, Santiago, Chile.
\email{mmatamal@dim.uchile.cl}
\and
Centro de Modelamiento Matem\'atico, UMR 2071, UCHILE-CNRS, Casilla
170-3, Correo 3, Santiago, Chile.
}

\maketitle

\begin{abstract}
Let be the following strategy to construct a walk in a labeled
digraph: at each vertex, we follow the unvisited arc of minimum label. 
In this work we study for which languages, applying the previous strategy
over the corresponding de Bruijn graph, we finish with an Eulerian
cycle, in order to obtain the minimal de 
Bruijn sequence of the language.
\end{abstract}

\section{Introduction}
Given a language, a de Bruijn sequence of span $n$  is a periodic
sequence such that every $n$-tuple in the language (and no other
$n$-tuple) occurs exactly once. Its first known description appears
as a Sanskrit word \textit{yam\'at\'ar\'ajabh\'anasalag\'am} which was
a memory aid for Indian drummers, where the accented/unaccented syllables
represent long/shorts beats, so all possible triplets of short and
long beats are included in the word. De Bruijn sequences are also known as
``shift register sequences'' and was originally studied by N. G. De
Bruijn for the binary alphabet \cite{deBruijn46:a_combinatorial}. These sequences have
many different applications, such as memory wheels in computers and
other technological device, network models, DNA algorithms,
pseudo-random number 
generation, modern public-key cryptographic schemes, to mention a few
(see \cite{SteinSciAmer},\cite{BermondBruijnKatz},\cite{ChungDiaconis}).
Historically, de Bruijn sequence was studied in an arbitrary alphabet considering the
language of all the $n$-tuples. There is a large number of de Bruijn
sequence in this case, but only a few can be generated efficiently,
see \cite{Fredricksen:82:SFL} for a survey about this subject. 
In 1978, Fredricksen and Maiorana \cite{Maiorana:77} give an
algorithm to generate a de Bruijn sequence of span $n$ based in the
Lyndon words of the language, which resulted to be the minimal one in
the lexicographic order, and this  algorithm was proved to be
efficient \cite{ruskey92:_gener_neckl}. 
Recently, the study of these concepts was extended to languages with
forbidden substrings: in \cite{Ruskey:2000:GNS} it was given efficient
algorithms to generate all the words in a language with one forbidden
substring, in \cite{moreno03:_lyndon_words_bruij_subsh_finit_type}
the concept of de Bruijn sequences was generalized to 
restricted languages with a finite set of forbidden substrings
and it was proved 
the existence of these sequences and presented an algorithm to generate
one of them, however, to find the minimal sequence is a non-trivial
problem in this more general case.  This problem is closely related to
the ``shortest common super-string problem''  which is a important
problem in the areas of DNA sequencing and data compression. 

In this work we study the de Bruijn sequence of minimal
lexicographical label.  
In section 2 we present some definitions and previous results on de
Bruijn sequences and the BEST Theorem, necessary to understand the
main problem, and we prove a result related with the BEST Theorem
which will be useful in the following sections. 
In section 3 we study the
main problem, giving some results on the structure of the de Bruijn
graph.
Finally, in section 4 we present some remarks and extensions
of this work.

\section{De Bruijn Sequence of Restricted Languages}

\subsection{Definitions}
Let $A$ be a finite set with a linear order $<$. A {\it word} on the
alphabet $A$ is a finite sequence of elements of $A$, whose length
is denoted by $|w|$. 

A word $p$ is said to be a {\it factor} of
a word $w$ if there exist words $u,v\in A^*$ such that $w=u p v$. 
If $u$ is the empty word  $\varepsilon$ then $p$ is
called a {\it prefix} of $w$, and if $v$ is empty then is called a
{\it suffix} of $w$. If $p\neq w$ then $p$ is a {\it proper factor},
{\it proper prefix} or {\it proper suffix}, respectively.

The set $A^*$ of all the words on the alphabet $A$ is linearly ordered by
the alphabetic order induced by the order $<$ on $A$. By definition,
$x<y$ either if $x$ is a prefix of $y$ or if $x=u a v$, $y=u b w$
with $u,v,w\in A^*$, $a,b\in A$ and $a<b$. A basic property of the
alphabetic order is the following: if $x<y$ and if $x$ is not a prefix of $y$,
then for any pair of words $u,v$, $x u<y v$.

Given an alphabet $A$, a full shift $A^\Zset$ is the collection of all
bi-infinite sequences of symbols from $A$.  Let $\mathcal{F}$ be a set
of words over $A^*$. A {\it subshift of finite type} (SFT) is the
subset of sequences in $A^\Zset$ which does not contain any factor in
$\mathcal{F}$.  We will refer to $\mathcal{F}$ as the set of {\it 
forbidden blocks} or {\it forbidden factors}.

Given a set $\mathcal{F}$ of forbidden blocks, in this work we will
say that a word $w$ is in the
language if the periodical word $w^\infty$, 
composed by infinite repetitions of $w$, is in the language of the SFT
defined by $\mathcal{F}$. The set of all the words of length $n$ in
the language defined by $\mathcal{F}$ will be 
denoted by $\conju{W}{n}$.

A SFT is {\it irreducible} if for every ordered pair of blocks $u,v$
in the language there is a block $w$ in the language  so that $u w v$
is a block of the language.

A de Bruijn sequence of span $n$  in a restricted language is a
circular string 
$\bru$ of length $\left|\conju{W}{n}\right|$ such that all the words
in the language of length $n$ are factors of $\bru$. In other words, 
\[\{ (\bru)_i\ldots (\bru)_{i+n-1 \mod n} | i=0\ldots n-1\}=\conju{W}{n}\] 

These concepts are studied in
\cite{moreno03:_lyndon_words_bruij_subsh_finit_type}, extending the
known results on subshifts of finite type to this context. In
particular  two results are relevant in this work, the first one is a
bound in the 
number of words of length $n$ in the language:
\[ \left|\conju{W}{n}\right| = \Theta\left(\lambda^n\right) \]
where $\log(\lambda)$ is the \textit{entropy} of the system (see
\cite{Lind:SDC}). The second result proves the existence of a de
Bruijn sequence:

\begin{theorem}\label{T:existenciaB}
For any set of forbidden substrings $\mathcal{F}$ defining an
irreducible subshift of finite type, there exists a de Bruijn sequence
of span $n$.
\end{theorem}

This last theorem is a direct consequence of the fact that the
de Bruijn graph of span $n$ is an 
Eulerian graph. The \textit{de Bruijn graph} of span $n$, denoted by $\gbru{n}$, is the
largest strongly connected component of the directed graph with
$|A|^n$ vertices, 
labeled by the words in  $A^n$, and the set of arcs  
\[ E=\left\{ (a s,s b)|a,b\in A, s\in A^{n-1}, a s b\in
  \conju{W}{n+1}\right\} \]
where the label of the arc $e=(a s , s b)$ is $l(e)=b$.
Note that if the SFT is irreducible, this graph has only one strongly connected
component of size greater than 1, so there is no ambiguity in the
definition. 

There are not two vertices with the same label, hence from now 
we identify a vertex by its label. If $W=e_1\ldots e_k$ is a walk over
$\gbru{n}$, we denote the label of $W$ by $l(W)=l(e_1)\ldots l(e_k)$,
and by $l(W)^j$ the concatenation of of $j$ times $l(W)$. 

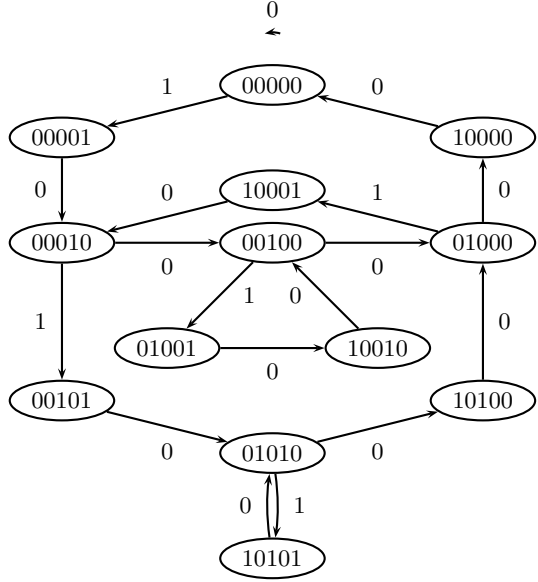
\begin{figure}[t]
\begin{center}
\begin{pspicture}(5.6,8)
\psset{unit=0.7cm}
\psset{arrows=->}
\rput(4,9){\ovalnode{0}{00000}}
\rput(0,8){\ovalnode{1}{00001}}
\rput(8,8){\ovalnode{16}{10000}}
\rput(4,7){\ovalnode{17}{10001}}
\rput(0,6){\ovalnode{2}{00010}}
\rput(4,6){\ovalnode{4}{00100}}
\rput(8,6){\ovalnode{8}{01000}}
\rput(2,4){\ovalnode{9}{01001}}
\rput(6,4){\ovalnode{18}{10010}}
\rput(0,3){\ovalnode{5}{00101}}
\rput(8,3){\ovalnode{20}{10100}}
\rput(4,2){\ovalnode{10}{01010}}
\rput(4,0){\ovalnode{21}{10101}}

\nccircle{0}{.5}\nbput{0}
\ncline{0}{1}\taput{1}
\ncline{1}{2}\tlput{0}
\ncline{2}{4}\tbput{0}
\ncline{2}{5}\tlput{1}
\ncline{4}{8}\tbput{0}
\ncline{4}{9}\trput{1}
\ncline{5}{10}\tbput{0}
\ncline{8}{16}\trput{0}
\ncline{8}{17}\taput{1}
\ncline{9}{18}\tbput{0}
\ncline{10}{20}\tbput{0}
\ncarc{10}{21}\trput{1}
\ncline{16}{0}\taput{0}
\ncline{17}{2}\taput{0}
\ncline{18}{4}\tlput{0}
\ncline{20}{8}\trput{0}
\ncarc{21}{10}\tlput{0}
\end{pspicture}
\caption{De Bruijn digraph of span 5 for the Golden Mean ($\mathcal{F}=\{11\}$)} \label{figura}
\end{center}
\end{figure}

There exists a bijection between the arcs of $\gbru{n}$ and the words in
$\conju{W}{n+1}$, because to each arc with label $a\in A$ with tail at
 $w'\in A^{n}$ we can associate the word $w' a$ which
is, by definition, a word in $\conju{W}{n+1}$. Equally if $w'a$ is a word of
$\conju{W}{n+1}$, with $a\in A$, then there exists a vertex $w'$ 
and an arc with tail at this vertex with label $a$.

Furthermore, if a word $w$ is a label of a walk from $u$ to $v$ then $v$ is
a suffix of length $n$ of $u w$. In the same way, if
$w\in\conju{W}{n+1}$
then there is a cycle $C$ in $\gbru{n}$ with label $l(C)$ such that
$l(C)^{\frac{n+1}{|C|}}=w$.

With all these properties it is easy to see that a de Bruijn
sequence of span $n+1$ is exactly the label of an Eulerian cycle over
$\gbru{n}$. 

\subsection{The BEST Theorem}

BEST is an acronym of N. G. de Bruijn, T. van
Aardenne-Ehrenfest, C. A. B. Smith and W. T. Tutte, the BEST Theorem
(see \cite{tutte84:_graph_theor}) gives a correspondence between Eulerian
cycles in a  digraph and its rooted trees converging to the root
vertex. 

Let $r$ be a vertex of an Eulerian digraph $G=(V,E)$, a
spanning tree converging to the root $r$ is a spanning tree such that
there exists a directed path from each vertex to the root. 

Given an Eulerian cycle starting at the root of an Eulerian digraph,
if for every vertex of $G$ we take the last arc with tail at
this vertex in the cycle then we obtain a spanning tree converging to
the root. Conversely, given a spanning tree converging to
the root, a walk over $G$ starting at the root and using
the arc in the tree only if all the arcs with tail at
this vertex  has been used, is an Eulerian
cycle. A walk over the graph of this kind will be called a walk
``avoiding the tree''.

The BEST Theorem proves that for every different spanning tree we have
a different Eulerian cycle. Therefore it also allows us to calculate the
exact number of Eulerian cycles on a digraph, which is given by
\[ C_\mathcal{F}= M_T \cdot \prod_{i=1}^{|V|} (d^+(v_i)-1)!\]
where $M_T$ is the number of rooted spanning trees converging to a
given vertex. We bound the second term by
$((\bar{d}^+ -1)!)^{|V|}$ where $\bar{d}^+$ is the mean of
the outgoing degrees over all the vertices, so we have a lower bound
to the number of de Bruijn sequences
\[C_\mathcal{F} = \Omega\left(\lfloor\lambda-1\rfloor
!^{\lambda^{n-1}}\right)\]
in particular, for a system with $\lambda\geq3$ the number of the
Bruijn sequences of span $n$ is exponential in the number of words in
the language of length $n-1$. In the systems with $3>\lambda>1$ this
bound  is generally also true, 
because the underestimated term $M_T$ is generally exponential, for
example, in the system without restrictions of alphabet
$\{0,1\}$, this term is equal to $2^{2^{n-1}}$.

Now, we define formally a walk ``avoiding a subgraph''. 
Let $r$ be any vertex. For each vertex 
$v\neq r$ in $\gbru{n}$ let $e_v$ be any arc starting at $v$. 
Let $H$ be the spanning subgraph of $\gbru{n}$ with arc set 
$\{e_v:v\in V(\gbru{n})\setminus \{r\}\}$.

Is easy to see that $H$ is composed by cycles, subtrees converging to
a cycle, and one subtree converging to $r$. For a vertex not in a
cycle of $H$, we define $H_v$ as the directed subtree converging to
$v$ in $H$.

We define recursively a walk in $\gbru{n}$ which \emph{avoid} $H$. 
It starts at the root vertex $r$.
Let $v_0e_0\cdots v_i$ be the current walk.
If there is an unvisited arc $e_i=(v_i,v_{i+1})$ not in $H$ 
we extend the walk by $e_i v_{i+1}$. Otherwise we use the arc
$e_{v_i}$ in $H$.

We say that a walk over the graph \textit{exhausts} a vertex if
the walk use all the arc having the vertex as head or tail. 

The next
lemma studies in which order the vertices are exhausted in a walk
avoiding $H$

\begin{lemma}\label{L:ExhaustSubtree}
Let $W$ be a walk starting at vertex $r$ avoiding $H$, let $v$ be a
vertex and let $W v$ the subpath of 
$W$ starting at vertex $r$ and finishing when it exhausts the
vertex $v$. Then 
for each vertex $u$ in $H_v$, $u$ is exhausted in
$W v$.
\end{lemma}

\begin{proof}
By induction in the depth of the subtree with root $v$. If $v$ is a
leaf of $H$ then $H_v=\{v\}$. If $v$ is not a
leaf and $W v$ exhaust $v$, then $W v$ visit all arc $(v,w)\in E$,
and therefore all the arcs $(u,v)\in E$, 
applying induction hypothesis to all vertices $u$ such that $(u,v)\in
E$ we prove the result.
\qed
\end{proof}

\section{Minimal de Bruijn Sequence}

Let $m=m_1,\ldots m_{n}$ be  the vertex of
$\gbru{n}$ of maximum label in the lexicographic order. We are
interested in to obtain the Eulerian cycle of minimum label starting
at $m$. In  order to obtain this cycle, we define the following walk: 
Starting at $m$, at each vertex  we continue by the arc with
the lowest  label between the unvisited arcs with tail at this
vertex. A walk constructed by this way will be called a
\textit{minimal walk}. By definition, there is no walk with a
lexicographically lower label, except its subwalks. 
In this section we characterize  
when a minimal walk starting at $m$ is an Eulerian cycle, obtaining the
minimal de Bruijn sequence.

For each vertex $v$ let $e(v)$ be the arc with tail at the vertex
$v$ and with maximum label. 
Let $T$ be the spanning subgraph of
$\gbru{n}$ composed by the set of arcs $e(v)$, for $v\in V(\gbru{n})$,
$v\neq m$. 
The label of $e(v)$ will be denoted by $\gamma(v)$.

Is easy to see that a minimal walk is a walk avoiding $T$, hence we
can study a minimal walk analyzing the structure of $T$.

\begin{theorem}\label{t:initial}
A minimal walk is an Eulerian cycle if and only if $T$ is a tree.
\end{theorem}
\begin{proof}
A minimal walk $W$ exhaust $m$, if $T$ is a tree then by Lemma
\ref{L:ExhaustSubtree} all vertices of $T$ are exhausted by
$W$, hence $W$ is an Eulerian cycle. Conversely, if $W$ is an Eulerian
cycle, by 
the BEST Theorem the subgraph composed by the last arc visited at each 
vertex is a tree, but this subgraph is $T$, concluding that $T$ is a
tree.
\qed\end{proof}

In the unrestricted case (when $\conju{W}{n}=A^n$), 
the subgraph $T$ is a regular tree of depth $n$ where each non-leaf
vertex has $|A|$ sons, therefore the minimal walk is an Eulerian cycle.

In the restricted case, 
we do not obtain necessarily an Eulerian cycle, because $T$
is not necessarily a spanning tree converging to the
root due to the existence of cycles.  

We will study the structure of the graph $\gbru{n}$ and the subgraph
$T$, 
specially the cycles in $T$. The main theorem of this section
characterizes the label of cycles in $T$, allowing us to characterize
the languages where the minimal walk is an Eulerian cycle.

First of all, we will prove some properties of the de Bruijn
graph to understand the structure of the arcs and cycles in $T$. 

\begin{lemma}\label{l:walklabel}
Let $k\geq n+2$. Let $W=v_0e_0v_1e_1\cdots e_{k-1}v_k$ be a walk in $T$. 
Then $l(e_0)\leq l(e_{n+1})$.
\end{lemma}
\begin{proof}
Since $v_{n}=l(e_0)\cdots l(e_{n-1})$ we have that 
$l(e_1)\cdots l(e_{n-1})l(e_{n})l(e_0)\in \conju{W}{n+1}$. Hence there
exists an arc $(v_{n+1},u)$ with label $l(e_0)$, where
$v_{n+1}=l(e_1)\cdots l(e_{n-1})l(e_{n})$. 
By the definition of $T$, $l(e_0)\leq \gamma(v_{n+1})=l(e_{n+1})$. 
\qed\end{proof} 

\begin{corollary}\label{c:cycleslabel}
Let $C$ be a cycle in $T$.
Then $|C|$ divides $n+1$. 
Moreover for every vertex $u$ in $C$,
$u\gamma(u)=l(C)^{\frac{n+1}{|C|}}$.
\end{corollary}
\begin{proof}
Let consider the walk $W=v_0e_0\cdots e_{|C|-1}v_{|C|}=v_0e_0\cdots 
e_{(n+1)|C|-1}v_0e_0v_1$ 
as $n+1$ repetitions of the cycle $C$. From Lemma \ref{l:walklabel} we have  
$l(e_0)\leq l(e_{n+1})\leq l(e_{2(n+1)})\leq l(e_{(n+1)|C|})=l(e_0)$.
Since we can start the cycle in any vertex we conclude that
$l(e_i)=l(e_{(n+1)+i})$ for every $i=0,\ldots, |C|-1$. Hence
$|C|$ divides $n+1$. The second conclusion comes from 
the fact that the label of any walk of length at most $n$
ending in a vertex $u$ is a suffix of $u$.
\qed\end{proof}

Let $u\neq m$ be a vertex. Among all the words which are prefix of $m$ 
and suffix of $u$, let $g(u)$ be the longest one (notice that 
$g(u)$ could be the empty word $\varepsilon$ and $|g(u)|<n$).
Let $\alpha(u)=m_{|g(u)|+1}$ be the letter following the end
of $g(u)$ in $m$.
 
Notice that in the unrestricted case, $|g(u)|$ is the distance over the
graph from the vertex $u$ to $m$.
This function will be essential in the study of $T$. The next
lemma give us a bound over the label of the arcs in terms of the
function $g(\cdot)$.

\begin{lemma}\label{l:height}
For all pairs of adjacent vertices  $u$ and $v$,
$l(u v)\leq \alpha(u)$. Moreover, 
if $l(u v)<\alpha(u)$ then $g(v)=\varepsilon$ and 
if $l(u v)=\alpha(u)$ then $g(v)=g(u)l(u v)$.
\end{lemma}

\begin{proof}
$g(u)$ is a suffix of $u$, and $u l(u v)\in\conju{W}{n+1}$, so
$g(u)l(u v)$ is a prefix of a word in $\conju{W}{n+1}$. Since 
$m$ is the maximal word and $g(u)$ is a prefix of $m$ we get $l(u
  v)\leq \alpha(u)$. 

If $l(u v)=\alpha(u)$ then $g(u)l(u v)$ is a prefix of $m$ and a suffix
of $v$. Hence $g(u)l(u v)$ is a suffix of $g(v)$.
Since by removing the last letter of a suffix of $v$ we obtain a
suffix of $u$ we conclude $g(v)=g(u)l(u v)$.

We show that if $g(v)\neq \varepsilon$ then $\alpha(u)\geq l(u v)$.
Let $g(v)=g'(v)l(u v)$, then $g'(v)$ is a suffix of $u$ and a prefix of $m$.
Hence $g'(v)$ is a suffix of $g(u)$. Therefore $g'(v)\alpha(u)$ is a factor 
of $m$. By the definition of $g(v)$ and the maximality of $m$
$g(v)$ is greater or equal (lexicographically) than $g'(v)\alpha(u)$.
We conclude that $\alpha(u)\geq l(u v)$.
\qed\end{proof}

In the unrestricted case, where $T$ is a tree of depth $n$, all the
arcs not in $T$ go to a leaf. In the general case we can define an
analog to the leaves.

We say that a vertex $u$ is a \textit{floor} vertex 
if $g(u)=\varepsilon$. Notice that in the unrestricted case the leaves
of $T$ are the floor vertices. We say that a vertex $u$ is a
\textit{restricted} vertex  
if $\gamma(u)<\alpha(u)$.

\begin{corollary}
If a cycle in $T$ contains $l$ \textit{restricted} vertices, then it
has exactly $l$ \textit{floor} vertices.
\end{corollary}
\begin{proof}
From Lemma \ref{l:height} we know that 
if a vertex $u$ is restricted then for every arc $(u,v)$ the vertex $v$ is a
floor vertex.
To conclude it is enough to see that  in $T$ an arc $(u,v)$ with $u$
unrestricted has label 
$\alpha(u)$. Then $v$ is not a floor vertex.
\qed\end{proof}

\begin{corollary}\label{c:path}
Let $P$ be a path in $T$ starting in a floor vertex, ending in a vertex $v$
and with unrestricted inner vertices. Then $l(P)=g(v)$.
\end{corollary}
\begin{proof}
We apply induction on the length of $P$. The case where the length of
$P$ is zero is direct since $v$ is a floor vertex.
Let us consider the case where $P$ has length at least 1. 
Since $v$ is not a restricted vertex, from Lemma \ref{l:height} we
know that  $g(v)=g(u)l(u v)$, where $u$ is
its neighbor in $P$. By the induction assumption $g(u)=l(P')$
where $P'$ is the path obtained from $P$ removing the arc $(u,v)$.
Hence $g(v)=l(P')l(u v)=l(P)$.
\qed\end{proof}

We will use these results to characterize the label of cycles in
$T$, specially we will characterize the restricted vertices of a cycle.

\begin{theorem}\label{L:SeqCycles}
Let $C$ be a cycle in $T$, let $u^0,\ldots, u^{k-1}$ be the \textit{restricted}
vertices in $C$ ordered according to the order of $C$.
Then $u^i=g(u^{i+1})\gamma(u^{i+1})\cdots \gamma(u^{i-1})g(u^i)$ for
$i=0,\ldots,k-1$, where $i+1,\ldots, i-1$ are computed $\mod k$.
\end{theorem} 
\begin{proof}
From Corollary \ref{c:path} the label of $C$ is
$g(u^0)\gamma(u^0)\cdots$ $g(u^{k-1})$ $\gamma(u^{k-1})$,
and by definition of $\gbru{n}$, $u^i$ is
the label of any walk over $\gbru{n}$ of length $n$ finishing in $u^i$, so by
Corollary  \ref{c:cycleslabel} we can take the walk $C^k$ composed by
$k=(n+1)/|C|$ repetitions of $C$ finishing in $u^i$, concluding that 
$u^i= g(u^{i+1}) \gamma(u^{i+1}) \cdots \gamma(u^l)$ $l(C^{k-1}) g(u_1)$
$\cdots$ $\gamma(u^{i-1})$ $g(u^i)$.
\qed\end{proof}

\begin{figure}[ht]
\begin{center}
\begin{pspicture}(0,-1.5)(9,5)
\psset{arrows=->}
\rput(4.5,4){\rnode{14}{1110}}
\rput(4.5,3){\rnode{7}{0111}}
\rput(2.5,2){\rnode{3}{0011}}
\rput(7.5,2){\rnode{11}{1011}}
\rput(1,1){\rnode{1}{0001}}
\rput(4,1){\rnode{9}{1001}}
\rput(6.5,1){\rnode{5}{0101}}
\rput(8.5,1){\rnode{13}{1101}}
\rput(0,0){\rnode{0}{0000}}
\rput(1.5,0){\rnode{8}{1000}}
\rput(3,0){\rnode{4}{0100}}
\rput(4.5,0){\rnode{12}{1100}}
\rput(6,0){\rnode{2}{0010}}
\rput(7.5,0){\rnode{10}{1010}}
\rput(9,0){\rnode{6}{0110}}
\ncline{0}{1}
\ncline{8}{1}
\ncline{4}{9}
\ncline{12}{9}
\ncline{2}{5}
\ncline{10}{5}
\ncline{6}{13}
\ncline{1}{3}
\ncline{9}{3}
\ncline{5}{11}
\ncline{3}{7}
\ncline{7}{14}
\ncline{13}{10}
\nccurve{11}{6}
\psset{linecolor=lightgray,angleA=200,angleB=320}
\nccircle[arrows=<-]{0}{.3}
\nccurve{8}{0}
\nccurve{4}{8}
\nccurve{12}{8}
\nccurve{2}{4}
\nccurve{10}{4}
\nccurve{6}{12}
\nccurve[angleA=0,angleB=40]{14}{12}
\nccurve[angleA=200,angleB=200,ncurv=2.4]{3}{6}
\nccurve[angleA=200,angleB=210,ncurv=1.8]{1}{2}
\ncline[angleA=180,angleB=270]{9}{2}
\nccurve[angleA=0,angleB=90]{5}{10}
\end{pspicture}
\end{center}
\caption{Example of the subgraph $T$ for $n=4$ and
  $\mathcal{F}=\{01111\}$ in a binary alphabet.}
\end{figure}
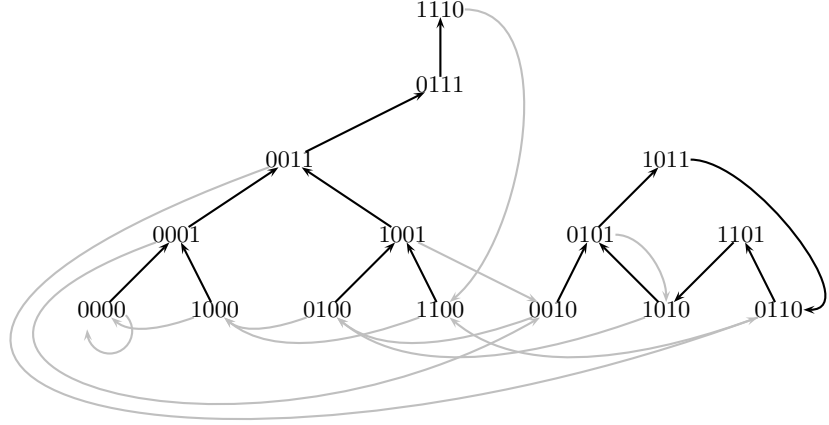

Now we are able to give a characterization of the languages where a
minimal walk produces an Eulerian cycle.

Let $\mathcal{H}$ be the subset of $\conju{W}{n+1}$ where $w\in
\mathcal{H}$ if and only if $w$ can be decomposed by
$w=h^0\beta_1\ldots h^{k-1}\beta_{k-1}$ where 
each $h^i\in A^*$ and $\beta_i\in A$ satisfy the following conditions:
\begin{enumerate}
\item $h^i=m_1\ldots m_{|h^i|}$ (a prefix of $m$)
\item $\beta_i<m_{|h^i|+1}$
\item $\forall \beta'>\beta_i,\  h^{i+1}\beta_{i+1}\ldots \beta_{i-1}h^i\beta'\notin \conju{W}{n+1}$ 
\end{enumerate}

Now, we are able to characterize the languages where a minimal walk is
an Eulerian cycle. 

\begin{theorem}
A minimal walk is an Eulerian cycle if and only if
$\mathcal{H}=\emptyset$.
\end{theorem}
\begin{proof}
From Theorem \ref{t:initial}, we have to prove that $T$ is a tree if
and only if $\mathcal{H}=\emptyset$.

If $T$ is not a tree then $T$ has a cycle $C$. 
Let  $u^0\ldots u^{k-1}$ be the restricted vertices of the cycle. 
By Theorem \ref{L:SeqCycles} $l(C)=g(u^0) \gamma(u^0) \ldots g(u^{k-1})
\gamma(u^{k-1})$ and 
by Corollary \ref{c:cycleslabel}
$|C|$ divides $n+1$. Therefore there exists a word
$w$ in $\conju{W}{n+1}$ composed by $(n+1)/|C|$ repetitions of $C$.
By definition of $\mathcal{H}$ we conclude that $w\in\mathcal{H}$.

Conversely, let us assume that $T$ has no cycles and $\mathcal{H}\neq
\emptyset$. Let $w$ be a word in $\mathcal{H}$. 
By definition of $\gbru{n}$, there is a cycle $C$ in $\gbru{n}$ of length
dividing $n+1$ such that $C$ (or repetitions of $C$) has label $w$.
We shall prove that $C$ is also a cycle in $T$. 

Let $v$ be a vertex of
$C$, with $v=\ldots \beta_{i-1}(h^i)_1\ldots (h^i)_j$ where
$j=0\ldots |h^i|$.
If $0<j<|h^i|$, then $m_1\ldots m_j$ is a suffix of $v$, so
$\alpha(v)=m_{j+1}=(h^i)_{j+1}$ hence the arc of $C$ with tail at $v$
is in $T$. 
If $j=0$ then $\gamma(v)=m_1$ therefore the arc in $C$ is in
$T$. 
Finally, let consider the case $j=|h^i|$.If  $(v,v')$ is the arc in
$C$ then $l(v v')=\beta_{i}$.
Since $w\in \mathcal{H}$, no arc in $\gbru{n}$ 
with tail at $v$ has a label greater than $\beta_i$. Then $(v,v')\in
T$.
We conclude that $C$ is a cycle in $T$ which leads to a
contradiction. 
\qed\end{proof}

\section{Some Remarks}

The previous analysis considers only the minimal walk starting at the root
vertex. This case does not necessarily produce the minimal label
over all Eulerian cycles, because there can be  Eulerian cycles
starting at a non root vertex with a lexicographically lower label. 

It is also possible to construct an algorithm which modifies $T$ in
order to destroy cycles in $T$, and obtain the minimal de Bruijn
sequence for any irreducible subshift of
finite type. However further research in this subject
allow us to construct an algorithm to obtain the minimal Eulerian
cycle for any edge-labeled digraph (see
\cite{matamala04:_minim_euler}), but this result escapes to the scope
of this work.

\bibliographystyle{splncs}
\bibliography{DinSimb,../CV/yo}

\begin{thebibliography}{10}

\bibitem{deBruijn46:a_combinatorial}
de~Bruijn, N.G.:
\newblock A combinatorial problem.
\newblock Nederl. Akad. Wetensch., Proc. \textbf{49} (1946)  758--764

\bibitem{SteinSciAmer}
Stein, S.K.:
\newblock The mathematician as an explorer.
\newblock Sci. Amer. \textbf{204} (1961)  148--158

\bibitem{BermondBruijnKatz}
Bermond, J.C., Dawes, R.W., Ergincan, F.{\"O}.:
\newblock De {B}ruijn and {K}autz bus networks.
\newblock Networks \textbf{30} (1997)  205--218

\bibitem{ChungDiaconis}
Chung, F., Diaconis, P., Graham, R.:
\newblock Universal cycles for combinatorial structures.
\newblock Discrete Math. \textbf{110} (1992)  43--59

\bibitem{Fredricksen:82:SFL}
Fredricksen, H.:
\newblock A survey of full length nonlinear shift register cycle algorithms.
\newblock SIAM Rev. \textbf{24} (1982)  195--221

\bibitem{Maiorana:77}
Fredricksen, H., Maiorana, J.:
\newblock Necklaces of beads in $k$ colors and $k$-ary de {B}ruijn sequences.
\newblock Discrete Math. \textbf{23} (1978)  207--210

\bibitem{ruskey92:_gener_neckl}
Ruskey, F., Savage, C., Wang, T.M.:
\newblock Generating necklaces.
\newblock J. Algorithms \textbf{13} (1992)  414--430

\bibitem{Ruskey:2000:GNS}
Ruskey, F., Sawada, J.:
\newblock Generating necklaces and strings with forbidden substrings.
\newblock Lect. Notes Comput. Sci. \textbf{1858} (2000)  330--339

\bibitem{moreno03:_lyndon_words_bruij_subsh_finit_type}
Moreno, E.:
\newblock Lyndon words and de bruijn sequences in a subshift of finite type.
\newblock In Harju, T., Karhum\"aki, J., eds.: Proceedings of WORDS'03.
  Number~27 in TUCS General Publications, Turku, Finland, Turku Centre for
  Computer Science (2003)  400--410

\bibitem{Lind:SDC}
Lind, D., Marcus, B.:
\newblock Symbolic Dynamics and Codings.
\newblock Cambridge University Press (1995)

\bibitem{tutte84:_graph_theor}
Tutte, W.T.:
\newblock Graph theory. Volume~21 of Encyclopedia of Mathematics and its
  Applications.
\newblock Addison-Wesley Publishing Company Advanced Book Program, Reading, MA
  (1984)

\bibitem{matamala04:_minim_euler}
Matamala, M., Moreno, E.:
\newblock Minimal {E}ulerian cycle in a labeled digraph.
\newblock Technical Report CMM-B-04/08-108, DIM-CMM, Universidad de Chile
  (2004)

\end{thebibliography}

\end{document}